# EMPIRICAL LIKELIHOOD FOR ESTIMATING EQUATIONS WITH MISSING VALUES

By Dong Wang and Song Xi Chen[1]

*University of Nebraska-Lincoln and Iowa State University and Peking University*

We consider an empirical likelihood inference for parameters defined by general estimating equations when some components of the random observations are subject to missingness. As the nature of the estimating equations is wide-ranging, we propose a nonparametric imputation of the missing values from a kernel estimator of the conditional distribution of the missing variable given the always observable variable. The empirical likelihood is used to construct a profile likelihood for the parameter of interest. We demonstrate that the proposed nonparametric imputation can remove the selection bias in the missingness and the empirical likelihood leads to more efficient parameter estimation. The proposed method is further evaluated by simulation and an empirical study on a genetic dataset on recombinant inbred mice.

**1. Introduction.** Missing data are encountered in many statistical applications. A major undertaking in biological research is to integrate data generated by different experiments and technologies. Examples include the effort by *genenetwork.org* and other data depositories to combine genetics, microarray data and phenotypes in the study of recombinant inbred mouse lines [34]. One problem in using measurements from multiple experiments is that different research projects choose to perform experiments on different subsets of mouse strains. As a result, only a portion of the strains have all the measurements, while other strains have missing measurements. The current practice of using only those complete measurements and ignoring incomplete measurements with missing values is undesirable since the selection bias in the missingness can cause the parameter estimators to be

Received November 2006; revised December 2007.

[1]Supported by NSF Grants SES-05-18904 and DMS-06-04563.

*AMS 2000 subject classifications.* Primary 62G05; secondary 62G20.

*Key words and phrases.* Empirical likelihood, estimating equations, Kernel estimation, missing values, nonparametric imputation.

inconsistent. Even in the absence of the selection bias (missing completely at random), the complete measurements-based inference is generally not efficient as it throws away data with missing values. Substantial research has been done to deal with missing data problems; see [17] for a comprehensive overview. Inference based on estimating equations [3, 9] is a general framework for statistical inference, accommodating a wide range of data structure and parameters. It has been used extensively for conducting semiparametric inference in the context of missing values. Robins, Rotnitzky and Zhao [25, 26] proposed using the parametrically estimated propensity scores to weigh estimating equations that define a regression parameter, and Robins and Rotnitzky [24] established the semiparametric efficiency bound for parameter estimation. The approach based on the general estimating equations has the advantage of being more robust against model misspecification, although a correct model for the conditional distribution of the missing variable given the observed variable is needed to attain the semiparametric efficiency. See [32] for a comprehensive review.

In this paper we consider an empirical likelihood based inference for parameters defined by general estimating equations in the presence of missing values. Empirical likelihood introduced by Owen [19, 20] is a computer-intensive statistical method that facilitates a likelihood-type inference in a nonparametric or semiparametric setting. It is closely connected to the bootstrap as the empirical likelihood effectively carries out the resampling implicitly. On certain aspects of inference, empirical likelihood is more attractive than the bootstrap, for instance its ability of internal studentizing so as to avoid explicit variance estimation and producing confidence regions with natural shape and orientation; see [21] for an overview. In an important development, Qin and Lawless [23] proposed an empirical likelihood for parameters defined by general estimating equations and established the Wilks theorem for the empirical likelihood ratio. Chen and Cui [5] show that the empirical likelihood of [23] is Bartlett correctable, indicating that the empirical likelihood has this delicate second-order property of the conventional likelihood more generally than previously anticipated. In the context of missing responses, Wang and Rao [33] studied empirical likelihood for the mean with imputed missing values from a kernel estimator of the conditional mean, and demonstrated that some of the attractive features of the empirical likelihood continue to hold.

When the parameter of interest defined by the general estimating equations is not directly related to a mean, or a regression model is not assumed as the model structure, the commonly used conditional mean-based imputation via either a parametric [36] or nonparametric [7] regression estimator may result in either biased estimation or reduced efficiency; for instance, when the parameter of interest is a quantile (conditional or unconditional) or some covariates are subject to missingness. To suit the general nature of parameters defined by general estimating equations and to facilitate a nonparametric likelihood inference in the presence of missing values, we propose a nonparametric imputation procedure that imputes missing values repeatedly from a kernel estimator of the conditional distribution of the missing variables given the fully observable variables. To control the variance of the estimating functions with imputed values, the estimating functions are averaged based on multiply-imputed values. We show that the maximum empirical likelihood estimator based on the nonparametric imputation is consistent and is more efficient than the estimator that ignores missing values. In particular, when the number of the estimating equations is the same as the dimension of the parameter, the proposed empirical likelihood estimator attains the semiparametric efficiency bound.

The paper is structured as follows. The proposed nonparametric imputation method is described in Section 2. The formulation of the empirical likelihood is outlined in Section 3. Section 4 gives theoretical results of the proposed nonparametric imputation-based empirical likelihood estimator. Results from simulation studies are reported in Section 5. Section 6 analyzes a genetic dataset on recombinant inbred mice. All technical details are provided in the Appendix.

**2. Nonparametric imputation.** Let $Z_i = (X_i^\tau, Y_i^\tau)^\tau$, $i = 1, \ldots, n$, be a set of independent and identically distributed random vectors, where $X_i$'s are $d_x$-dimensional and are always observable, and $Y_i$'s are $d_y$-dimensional and are subject to missingness. In practice, the missing components may vary among incomplete observations. For ease of presentation, we assume the missing components occupy the same components of $Z_i$. Extensions to the general case can be readily made. Furthermore, our use of $Y_i$ for the missing variable does not prevent it being either a response or covariates in a regression setting.

Let $\theta$ be a $p$-dimensional parameter so that $E\{g(Z_i, \theta)\} = 0$ at a unique $\theta_0$, which is the true parameter value. Here $g(Z, \theta) = (g_1(Z, \theta), \ldots, g_r(Z, \theta))^\tau$ represents $r$ estimating functions for an integer $r \geq p$. The interest of this paper is in the inference on $\theta$ when some $Y_i$'s are missing.

Define $\delta_i = 1$ if $Y_i$ is observed and $\delta_i = 0$ if $Y_i$ is missing. Like in [7, 33] and others, we assume that $\delta$ and $Y$ are conditionally independent given $X$, namely the strongly ignorable missing at random proposed by Rosenbaum and Rubin [27]. As a result,

$$P(\delta = 1 \mid Y, X) = P(\delta = 1 \mid X) =: p(X),$$

where $p(x)$ is the propensity score and prescribes a pattern of selection bias in the missingness.

Let $F(y|X_i)$ be the conditional distribution of $Y$ given $X = X_i$, and $W(\cdot)$ be a $d_x$-dimensional kernel function of the $q$th order satisfying

$$\int W(s_1, \ldots, s_{d_x}) \, ds_1 \cdots ds_{d_x} = 1,$$

$$\int s_i^l W(s_1, \ldots, s_{d_x}) \, ds_1 \cdots ds_{d_x} = 0 \qquad \text{for any } i = 1, \ldots, d_x \text{ and } 1 \leq l < q$$

and $\int s_i^q W(s_1, \ldots, s_{d_x}) \, ds_1 \cdots ds_{d_x} \neq 0$. A kernel estimator of $F(y|X_i)$ based on the sample is

$$(1) \qquad \hat{F}(y|X_i) = \sum_{l=1}^n \frac{\delta_l W((X_l - X_i)/h) I(Y_i \leq y)}{\sum_{j=1}^n \delta_j W((X_j - X_i)/h)}.$$

Here $h$ is a smoothing bandwidth and $I(\cdot)$ is the $d_y$-dimensional indicator function, which is defined as $I(Y_i \leq y) = 1$ if all components of $Y_i$ are less than or equal to the corresponding components of $y$, respectively, and $I(Y_i \leq y) = 0$ otherwise. The property of the kernel estimator when there are no missing values is well understood in the literature, for instance in [12]. Its properties in the context of the missing values can be established in a standard fashion. An important property that mirrors one for unconditional multivariate distribution estimators given in [15] is that the efficiency of $\hat{F}(y|X_i)$ is not influenced by the dimension of $Y_i$. Here we concentrate on the case that $X_i$ is a continuous random vector. Extension to discrete random variables can be readily made; see Section 5 for an implementation with binary random variables.

We propose to impute a missing $Y_i$ with a $\tilde{Y}_i$, which is randomly generated from the estimated conditional distribution $\hat{F}(y|X_i)$. Effectively $\tilde{Y}_i$ has a discrete distribution where the probability of selecting a $Y_l$ with $\delta_l = 1$ is

$$(2) \qquad \frac{W\{(X_l - X_i)/h\}}{\sum_{j=1}^n \delta_j W\{(X_j - X_i)/h\}}.$$

To control the variability of the estimating functions with imputed values, we make $\kappa$ independent imputations $\{\tilde{Y}_{i\nu}\}_{\nu=1}^\kappa$ from $\hat{F}(y|X_i)$ and use

$$(3) \qquad \tilde{g}_i(\theta) = \delta_i g(X_i, Y_i, \theta) + (1 - \delta_i) \kappa^{-1} \sum_{\nu=1}^\kappa g(X_i, \tilde{Y}_{i\nu}, \theta)$$

as the estimating function for the $i$th observation. Like the conventional multiple-imputation procedure of Rubin [28], to attain the best efficiency, $\kappa$ is required to converge to $\infty$. Our numerical experience indicates that setting $\kappa = 20$ worked quite well in our simulation experiments reported in Section 5. A theoretical justification for this choice can be drawn from a remark to Theorem 2 in the next section.

The way missing values are imputed depends critically on the nature of the parameter $\theta$ and the underlying statistical model. A popular imputation method is to impute a missing $Y_i$ by the conditional mean of $Y$ given $X = X_i$ as proposed in [36] under a parametric regression model and in [7] and [33] via the kernel estimator for the conditional mean. However, this conditional mean imputation may not work for a general parameter and a general model structure other than the regression model; for instance, when the parameter is a correlation coefficient or a conditional or unconditional quantile [1] where the estimating equation is based on a kernel smoothed distribution function. Nor is it generally applicable to missing covariates in a regression context. In contrast, the proposed nonparametric imputation is applicable for any parameter defined by estimating equations, and the way we impute the missing values is independent of the estimating equations. The latter is specially attractive as this separation of the imputation and the inference steps is considered a major advantage of the multiple imputation approach proposed by Rubin [28].

It should be noted that, when $\kappa \to \infty$, the proposed method is equivalent to imputing the estimating functions with missing $Y_i$'s using the Nadaraya–Watson estimator of $E\{g(X_i, Y_i, \theta)|X_i\}$,

$$(4) \qquad \hat{m}_g(X_i, \theta) = \int g(X_i, y, \theta) \, d\hat{F}(y|X_i).$$

The imputation of the estimating function has the imputation and inference steps intertwined together except in some special cases, for instance when $\theta$ is the mean of $Y_i$ as considered in Cheng [7] and Wang and Rao [33]. In that case, $g(Z, \theta) = (Y - \theta)$ and $\hat{m}_g(X_i, \theta)$ is a simple difference between the kernel estimator of $E(Y|X_i)$ and $\theta$, which effectively separates the imputation and the inference step. However, for a general estimation equation, the imputation and inference steps may not be separable. This means that, as the search for an estimator of $\theta$ is made through the parameter space, the imputation has to be repeated whenever there is a change in the $\theta$ value. This computational burden would be particularly severe for the empirical likelihood, and more so when a resampling-based approach, for instance the bootstrap, is used to profile the empirical likelihood ratio. In contrast, the proposed approach generates a fixed set of missing values. Once they are generated, the same algorithm for data without missing values can be used without reimputation.

The curse of dimension is an issue with the kernel estimator $\hat{F}(y|X_i)$. However, as demonstrated in Section 4, since the target of the inference is a finite-dimensional $\theta$ rather than the conditional distribution $F(y|X_i)$, the curse of dimension does not pose any leading order effects on $\theta$-estimation as long as the bias of the kernel estimator is controlled. When $d_x \geq 4$, controlling the bias requires the order of the kernel $q > 2$, the so-called high-order kernel, so that $\sqrt{n}h^q \to 0$ instead of $\sqrt{n}h^2 \to 0$ when a conventional second order kernel is used. Using a high-order kernel may occasionally cause $\hat{F}(y|X_i)$ not being a proper conditional distribution as the imputation probability weights in (2) may be negative in the tails. However, the occurrence of this phenomenon is rare for large sample sizes as $\hat{F}(y|x)$ is a consistent estimator of $F(y|x)$. In practice, we can readjust the probability weights in (2) by setting negative weights to zero and rescaling the remaining weights to assure that all weights sum up to one. This readjustment is similar to the method used by Hall and Murison [10] for high-order kernel density estimators.

## 3. Empirical likelihood.

The nonparametric imputation produces an extended sample including $(X_i, Y_i)^\tau$ for each $\delta_i = 1$, and $(X_i, \{\tilde{Y}_{i\nu}\}_{\nu=1}^\kappa)^\tau$ for each $\delta_i = 0$. With the imputed estimating functions $\tilde{g}_i(\theta)$, the usual estimating equation approach can be used to make inference on $\theta$. The variance of the general estimating equation-based estimator for $\theta$ can be estimated using a sandwich estimator and the confidence regions can be obtained by asymptotic normal approximation. In this article, we would like to carry out a likelihood type inference using empirical likelihood, encouraged by its attractive performance for estimating equations without missing values, as demonstrated by Qin and Lawless [23] and the work of Wang and Rao [33] for inference on a mean with missing responses. An advantage of empirical likelihood is that it has no predetermined shape of the confidence region; instead, it produces regions that reflect the features of the data set. Our proposal of using empirical likelihood in conjunction with the nonparametric imputation is especially attractive, since it requires only weak assumptions for both imputation and inference procedures while it also has the flexibility inherent to empirical likelihood and the estimating equations.

Let $p_i$ represent the probability weight allocated to $\tilde{g}_i(\theta)$. The empirical likelihood for $\theta$ based on $\tilde{g}_i(\theta)$ is

$$L_n(\theta) = \sup\left\{\prod_{i=1}^n p_i \middle| p_i \geq 0, \sum_{i=1}^n p_i = 1, \sum_{i=1}^n \tilde{g}_i(\theta) = 0\right\}.$$

By the standard derivation of empirical likelihood [23], the optimal $p_i$ is

$$p_i = \frac{1}{n}\frac{1}{1 + t^\tau(\theta)\tilde{g}_i(\theta)},$$

where $t(\theta)$ is the Lagrange multiplier that satisfies

$$(5) \qquad Q_{n1}(\theta) =: \frac{1}{n}\sum_i \frac{\tilde{g}_i(\theta)}{1 + t^\tau(\theta)\tilde{g}_i(\theta)} = 0.$$

Let $\ell_n(\theta) = -\log\{L_n(\theta)/n^{-n}\}$ be the log empirical likelihood ratio. The maximum empirical likelihood estimator (MELE), $\hat{\theta}_n$, can be derived by maximizing $L_n(\theta)$ or minimizing $\ell_n(\theta)$.

When the estimating function $g(Z, \theta)$ is differentiable with respect to $\theta$, the MELE can be found via solving the following system of equations [23],

$$(6) \qquad Q_{n1}(\theta, t) = 0 \quad \text{and} \quad Q_{n2}(\theta, t) = 0,$$

where $Q_{n1}(\theta, t)$ is given in (5) and

$$Q_{n2}(\theta, t) = \frac{1}{n} \sum_i \frac{1}{1 + t^\tau(\theta)\tilde{g}_i(\theta)} \left\{ \frac{\partial \tilde{g}_i(\theta)}{\partial \theta} \right\}^\tau t(\theta).$$

Like the conventional parametric maximum likelihood estimation (MLE), there may be multiple solutions to the likelihood equation (6) depending on the form of $g(Z, \theta)$ and the underlying distribution. It is required that each solution be substituted back to $L_n(\theta)$ to identify the MELE.

**4. Main results.** In this section, we first present a theorem regarding the consistency of $\hat{\theta}_n$, which is a solution of the likelihood equation (6). We then discuss the estimation efficiency and propose confidence regions for $\theta$ based on the empirical likelihood ratio. We use $\theta_0$ to denote the true parameter value.

THEOREM 1. *Under the conditions given in the Appendix, as $n \to \infty$ and $\kappa \to \infty$, with probability tending to 1 the likelihood equation (6) has a solution $\hat{\theta}_n$ within the open ball $\|\theta - \theta_0\| < Cn^{-1/3}$ for a positive constant.*

The theorem indicates consistency of $\hat{\theta}_n$. The nature of the result corresponds to Lemma 1 of Qin and Lawless [23] on the consistency of the maximum empirical likelihood estimator without missing values and is an analogue of Cramér [8] for parametric MLEs.

Next we consider the efficiency of $\hat{\theta}_n$. Write $g(Z) =: g(Z, \theta_0)$. We define

$$(7) \qquad \Gamma = E[p^{-1}(X) \operatorname{Var}\{g(Z)|X\} + E\{g(Z)|X\}E\{g^\tau(Z)|X\}],$$

$$(8) \qquad \tilde{\Gamma} = E[p(X) \operatorname{Var}\{g(Z)|X\} + E\{g(Z)|X\}E\{g^\tau(Z)|X\}]$$

and $V = \{E(\frac{\partial g}{\partial \theta})^\tau \tilde{\Gamma}^{-1} E(\frac{\partial g}{\partial \theta})\}^{-1}$ at $\theta = \theta_0$.

THEOREM 2. *Under the conditions given in the Appendix, as $n \to \infty$ and $\kappa \to \infty$,*

$$\sqrt{n}(\hat{\theta}_n - \theta_0) \xrightarrow{\mathcal{L}} N(0, \Sigma)$$

*with $\Sigma = VE(\frac{\partial g}{\partial \theta})^\tau \tilde{\Gamma}^{-1} \Gamma \tilde{\Gamma}^{-1} E(\frac{\partial g}{\partial \theta})V$.*

The estimator $\hat{\theta}_n$ is consistent and asymptotically normal for $\theta_0$ and the potential selection bias in the missingness as measured by the propensity score $p(x)$ has been filtered out. If there are no missing values, $\tilde{\Gamma} = \Gamma = E(gg^\tau)$, which means that

$$\Sigma = \left\{ E\left(\frac{\partial g}{\partial \theta}\right)^\tau (Egg^\tau)^{-1} E\left(\frac{\partial g}{\partial \theta}\right) \right\}^{-1}.$$

This is the asymptotic variance of the maximum empirical likelihood estimator based on full observations given in [23]. Comparing the forms of $\Sigma$ with and without missing values shows that the efficiency of the maximum empirical likelihood estimator based on the proposed imputation will be close to that based on full observations if either the proportion of missing data is low, or if $E\{p^{-1}(X)\operatorname{Var}(g|X)\}$ is small relative to $E\{E(g|X)E(g^\tau|X)\}$, namely when $X$ is highly "correlated" with $Y$.

In the case of $\theta = EY$, $\Sigma = E\{\sigma^2(X)/p(X)\} + \operatorname{Var}\{m(X)\}$, where $\sigma^2(X) = \operatorname{Var}(Y|X)$ and $m(X) = E(Y|X)$. Thus, in this case $\hat{\theta}_n$ is asymptotically equivalent to the estimator proposed by Cheng [7] and Wang and Rao [33] based on the conditional mean imputation.

When $r = p$, namely the number of estimating equations is the same as the dimension of $\theta$,

$$\Sigma = \left\{ E\left(\frac{\partial g}{\partial \theta}\right)^\tau \Gamma^{-1} E\left(\frac{\partial g}{\partial \theta}\right) \right\}^{-1},$$

which is the semiparametric efficiency bound for the estimation of $\theta$ as given by Chen, Hong and Tarozzi [6].

Like the multiple imputation of Rubin [28], our method requires $\kappa \to \infty$. To appreciate this proposal, we note that when $\kappa$ is fixed, the $\Gamma$ and $\tilde{\Gamma}$ matrices used to define $\Sigma$ are

$$\Gamma = E[\{p^{-1}(X) + \kappa^{-1}(1 - p(X))\}\operatorname{Var}(g|X) + E(g|X)E(g^\tau|X)]$$

and

$$\tilde{\Gamma} = E[\{p(X) + \kappa^{-1}(1 - p(X))\}\operatorname{Var}(g|X) + E(g|X)E(g^\tau|X)].$$

Hence, a larger $\kappa$ will reduce the terms in $\Gamma$ and $\tilde{\Gamma}$, which are due to a single nonparametric imputation. Our numerical experience suggests that $\kappa = 20$ is sufficient for most situations and would make the $\kappa^{-1}(1 - p(X))$-term small enough.

Let us now turn our attention to the log empirical likelihood ratio

$$\mathcal{R}_n(\theta_0) = 2\ell_n(\theta_0) - 2\ell_n(\hat{\theta}_n).$$

Let $I_r$ be the $r$-dimensional identity matrix. The next theorem shows that the log empirical likelihood ratio converges to a linear combination of independent chi-square distributions.

THEOREM 3. *Under the conditions given in the [Appendix], as $n \to \infty$ and $\kappa \to \infty$,*

$$\mathcal{R}_n(\theta_0) \xrightarrow{\mathcal{L}} Q^\tau \Omega Q,$$

*where $Q \sim N(0, I_r)$ and $\Omega = \Gamma^{1/2} \tilde{\Gamma}^{-1} E(\frac{\partial g}{\partial \theta}) V E(\frac{\partial g}{\partial \theta})^\tau \tilde{\Gamma}^{-1} \Gamma^{1/2}$.*

When there are no missing values, $\Gamma = \tilde{\Gamma} = E(gg^\tau)$ and

$$\Omega = E(gg^\tau)^{-1/2} E\left(\frac{\partial g}{\partial \theta}\right) \left[ E\left(\frac{\partial g}{\partial \theta}\right)^\tau \{E(gg^\tau)\}^{-1} E\left(\frac{\partial g}{\partial \theta}\right) \right]^{-1}$$
$$\times E\left(\frac{\partial g}{\partial \theta}\right)^\tau E(gg^\tau)^{-1/2},$$

which is symmetric and idempotent with $\mathrm{tr}(\Omega) = p$. This means that

$$\mathcal{R}_n(\theta_0) \xrightarrow{\mathcal{L}} \chi_p^2,$$

which is the nonparametric version of Wilks theorem established in Qin and Lawless [23].

When there are missing values, Wilks theorem for empirical likelihood is no longer available due to a mismatch between the variance of the quantity $n^{-1/2} \sum_{i=1}^n \tilde{g}_i(\theta_0)$ and the probability limit of $n^{-1} \sum_{i=1}^n \tilde{g}_i(\theta_0) \tilde{g}_i^\tau(\theta_0)$. This phenomenon also appears when a nuisance parameter is replaced by a plugged-in estimator as revealed by Hjort, McKeague and Van Keilegom [13].

When $\theta = EY$, $\mathcal{R}_n(\theta_0) \xrightarrow{\mathcal{L}} \{V_1(\theta_0)/V_2(\theta_0)\}\chi_1^2$, where

$$V_1(\theta_0) = E\{\sigma^2(X)/p(X)\} + \mathrm{Var}\{m(X)\}$$

and $V_2(\theta_0) = E\{\sigma^2(X)p(X)\} + \mathrm{Var}\{m(X)\}$. This is the limiting distribution given in [33].

As confidence regions can be readily transformed to test statistics for testing a hypothesis regarding $\theta$, we shall focus on confidence regions. There are potentially several methods for constructing confidence regions for $\theta$. One is based on an estimation of the covariance matrix $\Sigma$ and the asymptotic normality given in Theorem [2]. Another method is to estimate the matrix $\Omega$ in Theorem [3] and then use Fourier inversion or a Monte Carlo method to simulate the distribution of the linear combinations of chi-squares. Despite the loss of Wilks theorem, confidence regions based on the empirical likelihood ratio $R_n(\theta)$ still enjoy the attractions of likelihood-based confidence regions in terms of having natural shape and orientation and respecting the range of $\theta$.

We propose the following bootstrap procedure to approximate the distribution of $R_n(\theta_0)$. Bootstrap for imputed survey data has been discussed in Shao and Sitter [30] in the context of ratio and regression imputations. We use the following bootstrap procedure in which the bootstrap data set is imputed in the same way as the original data set:

1. Draw a simple random sample $\boldsymbol{\chi}_n^*$ with replacement from the extended sample $\bar{\boldsymbol{\chi}}_n = \{(X_i, Y_i)^\tau$ for each $\delta_i = 1$ and $(X_i, \{\tilde{Y}_{i\nu}\}_{\nu=1}^\kappa)^\tau$ for each $\delta_i = 0; i = 1, \ldots, n\}$.

2. Let $\boldsymbol{\chi}_{nc}^*$ be the portion of $\boldsymbol{\chi}_n^*$ without imputed values and $\boldsymbol{\chi}_{nm}^*$ be the set of vectors in the bootstrap sample with imputed values. Then replace all the imputed $Y$ values in $\boldsymbol{\chi}_{nm}^*$ using the proposed imputation method where the estimation of the conditional distribution is based on $\boldsymbol{\chi}_{nc}^*$.

3. Let $\ell^*(\hat{\theta}_n)$ be the empirical likelihood ratio based on the reimputed data set $\boldsymbol{\chi}_n^*$, $\hat{\theta}_n^*$ be the corresponding maximum empirical likelihood estimator and $\mathcal{R}^*(\hat{\theta}_n) = 2\ell^*(\hat{\theta}_n) - 2\ell^*(\hat{\theta}_n^*)$.

4. Repeat the above steps $B$-times for a large integer $B$ and obtain $B$ bootstrap values $\{\mathcal{R}_b^*(\hat{\theta}_n)\}_{b=1}^B$.

Let $q_\alpha^\star$ be the $1 - \alpha$ sample quantile based on $\{\mathcal{R}_b^*(\hat{\theta}_n)\}_{b=1}^B$. Then, an empirical likelihood confidence region with nominal coverage level $1 - \alpha$ is $I_\alpha = \{\theta \mid R(\theta) \leq q_\alpha^\star\}$. The following theorem justifies that this confidence region has correct asymptotic coverage.

THEOREM 4. *Suppose the conditions given in the Appendix are satisfied and $Q \sim N(0, I_r)$. Then, the conditional distribution of $\mathcal{R}^*(\hat{\theta}_n)$ given the original sample $\boldsymbol{\chi}_n$ converges to the distribution of $Q^\tau \Omega Q$ in probability as $n \to \infty$ and $\kappa \to \infty$.*

**5. Simulation results.** We report results from two simulation studies in this section. In each study, the proposed empirical likelihood inference based on the proposed nonparametric imputation is compared with the empirical likelihood inference based on (1) the *complete observations only* by ignoring data with missing values and (2) the *full observations* since the missing values are known in a simulation. When there is a selection bias in the missingness, the complete observations-based estimator may not be consistent. The proposed imputation will remove the selection bias in the missingness and improve estimation efficiency due to utilizing more data information. Obtaining the full observations-based estimator allows us to gauge how far away the proposed imputation based estimator is from the ideal case.

We also compare the proposed method with a version of the inverse probability weighted generalized method of moments (IPW-GMM) described in [6], in which the estimating functions involving complete observations are inflated by nonparametrically estimated propensity scores. Based on the usual formulation of the generalized method of moments (GMM) [11], the weighted-GMM estimator for $\theta_0$ considered in our simulation is

$$\tilde{\theta}_n = \arg\min_{\theta} \left\{ \frac{1}{n_c} \sum_{i=1}^{n} \delta_i g(Z_i, \theta) \frac{1}{\hat{p}(X_i)} \right\}^{\tau} A_T \left\{ \frac{1}{n_c} \sum_{i=1}^{n} \delta_i g(Z_i, \theta) \frac{1}{\hat{p}(X_i)} \right\},$$

where $n_c$ is the number of complete observations, $A_T$ is a nonnegative definite weighting matrix and $\hat{p}(X_i)$ is a consistent estimator for $p(X_i)$. The difference between the weighted-GMM method we use and that of [6] is that we used a kernel-based estimator for $p(X_i)$, instead of the sieve estimator described in [6]. The bandwidth used to construct $\hat{p}(X_i)$ is obtained by the cross-validation method. The kernel function $W(\cdot)$ is taken to be the Gaussian and product Gaussian kernels, respectively, for the two simulation studies. Cross-validation method is also used to choose the smoothing bandwidth in the kernel estimator $\hat{F}(y|X)$ given in (1) for the proposed nonparametric imputation; see [4] for details. Simulation results not reported here confirm that our proposed method is not sensitive to the choice of bandwidth. To satisfy the requirement $\sqrt{n}h^2 \to 0$, we use half of the bandwidth produced by the cross-validation procedure. This is only a rule of thumb. Alternatively, we could use the bandwidth obtained from the cross-validation with a higher order kernel. That would prescribe a bandwidth satisfying the condition asymptotically.

5.1. *Correlation coefficient.* In the first simulation, the parameter $\theta$ is the correlation coefficient between two random variables $X$ and $Y$ where $X$ is always observed, but $Y$ is subject to missingness. We first generate bivariate random vector $(X_i, U_i)^{\tau}$ from a skewed bivariate $t$-distribution suggested in [2] with five degrees of freedom, mean $(0,0)^{\tau}$, shape parameter $(4,1)^{\tau}$ and dispersion matrix

$$\bar{\Omega} = \begin{bmatrix} 1 & 0.955 \\ 0.955 & 1 \end{bmatrix}.$$

Then we let $Y_i = U_i - 1.2X_i I(X_i < 0)$. These make $(X_i, Y_i)^{\tau}$ have mean $(0, 0.304)^{\tau}$ and correlation coefficient $0.676$.

We consider three missing mechanisms:
(a) $p(x) = (0.3 + 0.175|x|)I(|x| < 4) + I(|x| \geq 4)$;
(b) $p(x) \equiv 0.65$ for all $x$;
(c) $p(x) = 0.5I(x > 0) + I(x \leq 0)$.
The mechanism (b) is missing completely at random, whereas the other two are missing at random and prescribe selection bias in the missingness.

Let $\mu_x$ and $\mu_y$ be the means, and $\sigma_x^2$ and $\sigma_y^2$ be the variances of $X$ and $Y$, respectively. In the construction of the empirical likelihood for $\theta$, the correlation coefficient, $\lambda = (\mu_x, \mu_y, \sigma_x^2, \sigma_y^2)^{\tau}$ are treated as nuisance parameters. When all observations are complete (no missing data), the estimating equation can be written as $n^{-1} \sum_{i=1}^{n} g(X_i, Y_i, \theta, \lambda) = 0$ with

$$g(X_i, Y_i, \theta, \lambda) = \begin{pmatrix} X_i - \mu_x \\ Y_i - \mu_y \\ (X_i - \mu_x)^2 - \sigma_x^2 \\ (Y_i - \mu_y)^2 - \sigma_y^2 \\ (X_i - \mu_x)(Y_i - \mu_y) - \theta \sigma_x \sigma_y \end{pmatrix}.$$

Table 1 contains the bias and standard deviation of the four estimators considered based on 1000 simulations with the sample size $n = 100$ and 200, respectively. It also contains the empirical likelihood confidence intervals using the full observations, complete observations only and the proposed nonparametric imputation method at a nominal level of 95%. They are all based on the proposed bootstrap calibration method with $B = 1000$. When using the nonparametric imputation method, $\kappa = 20$ imputations were made for each missing $Y_i$. The confidence intervals based on the weighted-GMM are calibrated using the asymptotic normal approximation with the covariance matrix estimated by the kernel method.

The results in Table 1 can be summarized as follows. The proposed empirical likelihood estimator based on the nonparametric imputation method significantly reduced the bias compared to inference based only on complete observations when the data were missing at random but not missing completely at random. The estimator based on the completely observed data suffered severe bias under missing mechanisms (a) and (c). The proposed estimator had smaller standard deviations than the complete observation-based estimator under all three missing mechanisms, including the case of missing completely at random. The weighted-GMM method also performed better than the complete observation-based estimator. However, it had larger variance than the proposed estimator. Most strikingly, the standard deviations of the empirical likelihood estimator based on the proposed imputation method were all quite close to the full observation-based estimator, which confirmed its good theoretical properties. Confidence intervals based on the complete observations only and the weighted-GMM method could have severe under-coverage: the former is due to the selection bias and the latter is due to the normal approximation. The proposed confidence intervals had satisfactory coverages, which are quite close to the nominal level 0.95.

5.2. *Generalized linear models with missing covariates.* In the second simulation study we consider missing covariates in a generalized linear model (GLM). We also take the opportunity to discuss an extension of the proposed imputation procedure to binary random variables. Commonly used methods in dealing with missing data for GLM are reviewed in [14]. Empirical likelihood for GLMs with no missing data was first studied by Kolaczyk

TABLE 1

*Inference for the correlation coefficient with missing values. The four methods considered are empirical likelihood using full observations, empirical likelihood using only complete observations (Complete obs.), inverse probability weighting based generalized method of moments (Weighted-GMM), and empirical likelihood using the proposed nonparametric imputation (N. imputation). The nominal coverage probability of the confidence interval is 0.95*

| Methods | Bias | Std. dev. | MSE | Coverage | Length of CI |
|---|---|---|---|---|---|
| | | $n = 100$ | | | |
| Full observations | −0.0026 | 0.0895 | 0.0080 | 0.936 | 0.3555 |
| | | Missing mechanism (a) | | | |
| Complete obs. | 0.0562 | 0.1222 | 0.0181 | 0.851 | 0.4967 |
| Weighted-GMM | 0.0108 | 0.1112 | 0.0125 | 0.776 | 0.2495 |
| N. imputation | −0.0092 | 0.1041 | 0.0109 | 0.945 | 0.4875 |
| | | Missing mechanism (b) | | | |
| Complete obs. | −0.0080 | 0.1162 | 0.0136 | 0.930 | 0.4482 |
| Weighted-GMM | −0.0150 | 0.1069 | 0.0117 | 0.802 | 0.2763 |
| N. imputation | −0.0138 | 0.0999 | 0.0101 | 0.932 | 0.4173 |
| | | Missing mechanism (c) | | | |
| Complete obs. | −0.1085 | 0.1442 | 0.0326 | 0.832 | 0.5593 |
| Weighted-GMM | −0.0266 | 0.1167 | 0.0143 | 0.786 | 0.2860 |
| N. imputation | −0.0383 | 0.1053 | 0.0125 | 0.928 | 0.4322 |
| | | $n = 200$ | | | |
| Full observations | 0.0071 | 0.0610 | 0.0038 | 0.958 | 0.2484 |
| | | Missing mechanism (a) | | | |
| Complete obs. | 0.0710 | 0.0776 | 0.0111 | 0.824 | 0.3161 |
| Weighted-GMM | 0.0112 | 0.0734 | 0.0055 | 0.799 | 0.2060 |
| N. imputation | 0.0038 | 0.0709 | 0.0050 | 0.955 | 0.3180 |
| | | Missing mechanism (b) | | | |
| Complete obs. | −0.0030 | 0.0799 | 0.0064 | 0.937 | 0.3091 |
| Weighted-GMM | −0.0031 | 0.0719 | 0.0052 | 0.832 | 0.2075 |
| N. imputation | −0.0023 | 0.0668 | 0.0045 | 0.942 | 0.2797 |
| | | Missing mechanism (c) | | | |
| Complete obs. | −0.0915 | 0.0979 | 0.0179 | 0.788 | 0.3919 |
| Weighted-GMM | −0.0107 | 0.0745 | 0.0057 | 0.820 | 0.2131 |
| N. imputation | −0.0118 | 0.0680 | 0.0048 | 0.936 | 0.2860 |

[16]. Application of empirical likelihood method to GLMs can help overcome difficulties with parametric likelihood, especially in the aspect of overdispersion.

To demonstrate how to extend the proposed method to discrete component in $X_i$, we consider a logistic regression model with binary response variable $X_3$ and covariates $X_1$, $X_2$ and $Y$. We choose logit$\{P(X_{3i} = 1)\} =$ $-1 + X_{1i} + X_{2i} - 1.5Y_i$, $X_{1i} \sim N(0, 0.5^2)$, $X_{2i} \sim N(3, 0.5^2)$ and $Y_i$ being binary with $\text{logit}\{P(Y_i = 1)\} = -1 + X_{1i} + 0.5X_{2i}$. Here $X_{1i}$, $X_{2i}$ and $X_{3i}$ are always observable while the binary $Y_i$ is subject to missingness with $\text{logit}\{P(Y_i \text{ is missing})\} = 0.5 + 2X_{1i} + X_{2i} - 3X_{3i}$. This model with $d_x = 3$ also allows us to see if there is a presence of the curse of dimension due to the use of the kernel estimator in the proposed imputation procedure.

When no missing data are involved, the empirical likelihood analysis for the logistic model simply involves the estimating equations

$$n^{-1} \sum_{i=1}^{n} S_i\{X_{3i} - \pi(S_i^\tau \beta)\} = 0$$

with $S_i = (1, X_{1i}, X_{2i}, Y_i)^\tau$, $\beta$ being the parameter and $\pi(z) = \exp(z)/\{1 + \exp(z)\}$. Although our proposed imputation in Section 2 is formulated directly for continuous random variables, binary response $X_{3i}$ can be accommodated by splitting the data into two parts according to the value of $X_{3i}$, and then applying the proposed imputation scheme to each part by smoothing on the continuous $X_{1i}$ and $X_{2i}$. The maximum empirical likelihood estimator for $\beta$ uses a modified version of the fitting procedure described in Chapter 2 of [18].

The results of the simulation study with $n = 150$ and 250 are shown in Table 2. Despite that the dimension of $X_i$ is increased to 3, the standard deviations of the proposed estimator were still quite close to the full observation-based empirical likelihood estimator, which was very encouraging. For parameters $\beta_0$, $\beta_1$ and $\beta_2$, the mean squared error of the proposed estimator is several folds smaller than that based on the complete observations only; the proposed method also leads to a reduction in the mean squared error by as much as one fold relative to the weighted-GMM. All three methods give similar mean squared errors for the parameter $\beta_3$, while the proposed estimator had the smallest mean squared error. The confidence intervals based on only complete observations or the weighted-GMM tend to show notable undercoverage, while the proposed confidence intervals have satisfactory coverage levels for all parameters.

**6. Empirical study.** Microarray technology provides a powerful tool in molecular biology by measuring the expression level of thousands of genes simultaneously. One problem of interest is to test whether the expression level of genes is related to a traditional trait like body weight, food consumption or bone density. This is usually the first step in uncovering roles that a gene plays in important pathways. The BXD recombinant inbred strains of mice were derived from crosses between C57BL/6J (B6 or B) and DBA/2J (D2 or D) strains [35]. Around 100 BXD strains have been established by researchers at University of Tennessee and the Jackson Laboratory. A variety of phenotype data are accumulated for a BXD mouse over the years [22].

TABLE 2

*Inference for parameters in a logistic regression model with missing values. The four methods considered are empirical likelihood using full observations (Full obs.), empirical likelihood using only complete observations (Complete obs.), inverse probability weighting based generalized method of moments (Weighted-GMM), and empirical likelihood using the proposed nonparametric imputation (N. imputation). The nominal coverage probability of the confidence interval is 0.95*

| Methods | Bias | Std. dev. | MSE | Coverage | Length of CI |
|---|---|---|---|---|---|
| | | $n = 150$ | | | |
| | | $\beta_0 = -1$ | | | |
| Full obs. | $-0.0296$ | 1.292 | 1.669 | 0.964 | 5.477 |
| Complete obs. | $-1.715$ | 1.618 | 5.559 | 0.920 | 6.840 |
| Weighted-GMM | $-0.7835$ | 1.562 | 3.053 | 0.891 | 5.250 |
| N. imputation | 0.0349 | 1.317 | 1.736 | 0.967 | 5.549 |
| | | $\beta_1 = 1$ | | | |
| Full obs. | 0.0519 | 0.4384 | 0.1949 | 0.964 | 1.820 |
| Complete obs. | 0.7898 | 0.5603 | 0.9377 | 0.796 | 2.510 |
| Weighted-GMM | 0.4302 | 0.5486 | 0.4860 | 0.834 | 1.811 |
| N. imputation | $-0.0605$ | 0.4388 | 0.1962 | 0.961 | 1.851 |
| | | $\beta_2 = 1$ | | | |
| Full Obs. | 0.0367 | 0.4500 | 0.2038 | 0.972 | 2.007 |
| Complete obs. | 0.4205 | 0.5590 | 0.4892 | 0.945 | 2.599 |
| Weighted-GMM | 0.2542 | 0.5484 | 0.3653 | 0.896 | 1.791 |
| N. imputation | $-0.0110$ | 0.4576 | 0.2095 | 0.966 | 1.993 |
| | | $\beta_3 = -1.5$ | | | |
| Full obs. | $-0.0531$ | 0.4979 | 0.2507 | 0.976 | 2.137 |
| Complete obs. | $-0.0684$ | 0.5713 | 0.3310 | 0.975 | 2.592 |
| Weighted-GMM | $-0.0751$ | 0.5843 | 0.3471 | 0.838 | 1.574 |
| N. imputation | 0.0718 | 0.5521 | 0.3100 | 0.966 | 2.474 |

The trait that we consider is the fresh eye weight measured on 83 BXD strains by Zhai, Lu and Williams (ID 10799, BXD phenotype data base). The Hamilton Eye Institute Mouse Eye M430v2 RMA Data Set contains measures of expression in the eye on 39,000 transcripts. It is of interest to test whether the fresh eye weight is related to the expression level of certain genes. However, the microarray data are only available for 45 out of the 83 BXD mouse strains for which fresh eye weights are all available. The most common practice is to use only complete observations and ignore missing values in the statistical inference. As demonstrated in our simulation, this approach can lead to inconsistent parameter estimators if there is a selection bias in the missingness. Even in the absence of selection bias, the estimators are not efficient as only those complete observations are used.

TABLE 2
*(Continued)*

| Methods | Bias | Std. dev. | MSE | Coverage | Length of CI |
|---|---|---|---|---|---|
| | | $n = 250$ | | | |
| | | $\beta_0 = -1$ | | | |
| Full obs. | $-0.0286$ | 0.9651 | 0.9321 | 0.956 | 3.916 |
| Complete obs. | $-1.670$ | 1.212 | 4.255 | 0.801 | 4.790 |
| Weighted-GMM | $-0.6393$ | 1.150 | 1.7304 | 0.862 | 3.832 |
| N. imputation | 0.0284 | 0.9801 | 0.9615 | 0.962 | 3.963 |
| | | $\beta_1 = 1$ | | | |
| Full obs. | 0.0195 | 0.3332 | 0.1114 | 0.953 | 1.349 |
| Complete obs. | 0.7270 | 0.4398 | 0.7220 | 0.665 | 1.789 |
| Weighted-GMM | 0.3166 | 0.4223 | 0.2786 | 0.782 | 1.304 |
| N. imputation | $-0.0660$ | 0.3367 | 0.1177 | 0.947 | 1.380 |
| | | $\beta_2 = 1$ | | | |
| Full obs. | 0.0305 | 0.3374 | 0.1147 | 0.958 | 1.400 |
| Complete obs. | 0.3902 | 0.4134 | 0.3232 | 0.867 | 1.729 |
| Weighted-GMM | 0.1966 | 0.3993 | 0.1981 | 0.874 | 1.297 |
| N. imputation | $-0.0173$ | 0.3406 | 0.1163 | 0.967 | 1.384 |
| | | $\beta_3 = -1.5$ | | | |
| Full obs. | $-0.0611$ | 0.3818 | 0.1495 | 0.950 | 1.529 |
| Complete obs. | $-0.0351$ | 0.4445 | 0.1988 | 0.963 | 1.797 |
| Weighted-GMM | $-0.0419$ | 0.4596 | 0.2130 | 0.791 | 1.165 |
| N. imputation | 0.0762 | 0.4377 | 0.1974 | 0.944 | 1.759 |

We conduct four separate simple linear regression analyses of the eye weight $(x)$ on the expression level $(y)$ of four genes, respectively. The estimating equation can be written as $n^{-1} \sum_{i=1}^{n} g(X_i, Y_i, \theta) = 0$, where

$$g(X_i, Y_i, \theta) = \begin{pmatrix} X_i - \theta_1 - \theta_2 Y_i \\ X_i Y_i - \theta_1 Y_i - \theta_2 Y_i^2 \end{pmatrix}$$

and $\theta_1$ and $\theta_2$ represent the intercept and slope, respectively. The genes are *H3071E5*, *Slc26a8*, *Tex9* and *Rps16*. Here we have missing covariates in our analysis. The missing gene expression levels are imputed from a kernel estimator of the conditional distribution of the gene expression level given the fresh eye weight. The smoothing bandwidths were selected based on the cross-validation method, which is 1.5 for the first three genes in Table 3 and 1.8 for gene *Rps16*.

Table 3 reports empirical likelihood estimates of the intercept and slope parameters and their 95% confidence intervals based on the proposed non-parametric imputation and empirical likelihood. It also contains results from a conventional parametric regression analysis using only the complete observations, assuming independent and identically normally distributed resid-uals. Table 3 shows that these two inference methods can produce quite different parameter estimates and confidence intervals. The difference in parameter estimates is as large as 50% for the intercept and 25% for the slope parameter. Table 3 also reports estimates and confidence intervals of the correlation coefficients using the proposed method and Fisher's $z$ transformation. The latter is based on the complete observations only and is the method used by *genenetwork.org*. We observe again differences between the two methods despite not being significant at the 5% level. The largest difference of about 30% is registered at gene *H3071E5*. As indicated earlier, part of the differences may be the estimation bias of the complete observations-based estimators as they are unable to filter out selection bias in the missingness.

TABLE 3
*Parameter estimates and confidence intervals (shown in parentheses) based on a simple linear regression model using the parametric method with complete observations only and the empirical likelihood method using the proposed nonparametric imputation. For the parametric inference, the confidence intervals for the intercept and slope are obtained using quantiles of the t-distribution, and the confidence intervals for the correlation coefficient are obtained by Fisher's z transformation*

| Gene | Complete observations only (parametric) | | Nonparametric imputation (with empirical likelihood) | |
|---|---|---|---|---|
| | | Intercept | | |
| *H3071E5* | −21.99 | (−40.97, −2.998) | −15.69 | (−37.02, 5.209) |
| *Slc26a8* | 73.59 | (49.45, 97.73) | 67.28 | (38.34, 95.87) |
| *Tex9* | −23.81 | (−46.12, −1.507) | −14.66 | (−38.57, 8.776) |
| *Rps16* | −13.52 | (−31.08, 4.041) | −8.090 | (−26.76, 10.18) |
| | | Slope | | |
| *H3071E5* | 10.16 | (5.720, 14.59) | 8.736 | (2.688, 14.21) |
| *Slc26a8* | −6.352 | (−9.294, −3.411) | −5.561 | (−9.431, −1.471) |
| *Tex9* | 5.101 | (2.588, 7.613) | 4.094 | (0.8753, 6.979) |
| *Rps16* | 6.766 | (3.371, 10.16) | 5.754 | (1.948, 9.236) |
| | | Correlation coefficient | | |
| *H3071E5* | 0.5757 | (0.3395, 0.7436) | 0.4426 | (0.1321, 0.6977) |
| *Slc26a8* | −0.5533 | (−0.7285, −0.3102) | −0.4319 | (−0.6809, −0.0761) |
| *Tex9* | 0.5296 | (0.2996, 0.7124) | 0.4024 | (0.1013, 0.6846) |
| *Rps16* | 0.5256 | (0.2744, 0.7097) | 0.4151 | (0.0755, 0.6613) |

## APPENDIX

Let $f(x)$ be the probability density function of $X$ and

$$m_g(x) = E\{g(X, Y, \theta_0)|X = x\}.$$

The following conditions are needed in the proofs of the theorems.

C1: The missing propensity function $p(x)$, the $X$-density $f(x)$ and $m_g(x)$ all have bounded partial derivatives with respect to $x$ up to an order $q$ with $q \geq 2$ and $2q > d_x$, and $\inf_x p(x) \geq c_0$ for some $c_0 > 0$.

C2: The true parameter value $\theta_0$ is the unique solution to $E\{g(Z, \theta) = 0\}$; the estimating function $g(x, y, \theta_0)$ has bounded $q$th order partial derivatives with respect to $x$, and $E\|g(Z, \theta_0)\|^4 < \infty$.

C3: The second partial derivative $\partial^2 g(z, \theta)/\partial\theta\,\partial\theta^\tau$ is continuous in $\theta$ in a neighborhood of the true value $\theta_0$; functions $\|\partial g(z, \theta)/\partial\theta\|$, $\|g(z, \theta)\|^3$ and $\|\partial^2 g(z, \theta)/\partial\theta\,\partial\theta^\tau\|$ are all uniformly bounded by some integrable function $M(z)$ in the neighborhood of $\theta_0$.

C4: The matrices $\Gamma$ and $\tilde{\Gamma}$ defined in (7) and (8) are, respectively, positive definite and $E[\partial g(z, \theta)/\partial\theta]$ has full column rank $p$.

C5: The smoothing bandwidth $h$ satisfies $nh^{d_x} \to \infty$ and $\sqrt{n}h^q \to 0$ as $n \to \infty$. Here $q$ is the order of the kernel $K$.

Assuming $p(x)$ being bounded away from zero in C1 implies that data cannot be missing with probability 1 anywhere in the domain of the $X$ variable. Conditions C2, C3 and C4 are standard assumption for empirical likelihood-based inference with estimating equations. In condition C5, that $\sqrt{n}h^q \to 0$ is to control the bias induced by the kernel smoothing, whereas that $nh^{d_x} \to \infty$ leads to consistent estimation of the conditional distribution. To simplify the exposition and without loss of generality, we will only deal with the case that $q = 2$ in the proof.

LEMMA A.1. *Assume that conditions* C1–C5 *are satisfied, then as* $n \to \infty$ *and* $\kappa \to \infty$,

$$n^{-1/2} \sum_{i=1}^{n} \tilde{g}_i(\theta_0) \xrightarrow{\mathcal{L}} N(0, \Gamma),$$

*where* $\Gamma = E\{p^{-1}(X)\operatorname{Var}(g|X) + E(g|X)E(g^\tau|X)\}$.

For the proof of Lemma A.1, we need the following proposition, which is a direct consequence of Lemma 1 in [29].

PROPOSITION A.1. *Let* $\{V_i\}$ *be a sequence of random variables such that, for some function* $h$, *as* $n \to \infty$, $h(V_1, \ldots, V_n) \xrightarrow{\mathcal{L}} \Xi$, *where* $\Xi$ *has a distribution function* $G$. *If* $\{U_i\}$ *is a sequence of random variables such that*

$$P\{U_n - h(V_1, \ldots, V_n) \leq s \mid V_1, \ldots, V_n\} \to F(s)$$

*almost surely for all* $s \in \mathbb{R}$, *where* $F$ *is a continuous distribution function, then*

$$P(U_n \leq t) \to (G * F)(t)$$

*for all* $t \in \mathbb{R}$, *where "*$*$" *denotes the convolution operator.*

PROOF OF LEMMA A.1. Let $u \in \mathbb{R}^r$ and $\|u\| = 1$. Also let $g_u(Z, \theta_0) = u^\tau g(Z, \theta_0)$ and $\tilde{g}_{ui}(\theta_0) = u^\tau \tilde{g}_i(\theta_0)$. First we need to show that

$$n^{-1/2} \sum_{i=1}^n \tilde{g}_{ui}(\theta_0) \overset{\mathcal{L}}{\to} N(0, u^\tau \Gamma u)$$

and then use the Cramér–Wold device to prove Lemma A.1. Define

$$m_{g_u}(x) = E(g_u(X, Y, \theta_0) | X = x)$$

and

$$\hat{m}_{g_u}(x) = \frac{\sum_{i=1}^n \delta_i W((x - X_i)/h) g_u(x, Y_i, \theta_0)}{\sum_{i=1}^n \delta_i W((x - X_i)/h)}.$$

Now we have

$$\frac{1}{n} \sum_{i=1}^n \left\{ \delta_i g_u(X_i, Y_i, \theta_0) + (1 - \delta_i) \kappa^{-1} \sum_{\nu=1}^\kappa g_u(X_i, \tilde{Y}_{i\nu}, \theta_0) \right\}$$

$$= \frac{1}{n} \sum_{i=1}^n \delta_i \{ g_u(X_i, Y_i, \theta_0) - m_{g_u}(X_i) \}$$

$$+ \frac{1}{n} \sum_{i=1}^n (1 - \delta_i) \left\{ \kappa^{-1} \sum_{\nu=1}^\kappa g_u(X_i, \tilde{Y}_{i\nu}, \theta_0) - \hat{m}_{g_u}(X_i) \right\}$$

$$+ \frac{1}{n} \sum_{i=1}^n (1 - \delta_i) \{ \hat{m}_{g_u}(X_i) - m_{g_u}(X_i) \} + \frac{1}{n} \sum_{i=1}^n m_{g_u}(X_i)$$

$$:= S_n + A_n + T_n + R_n.$$

Note that $S_n$ and $R_n$ are sums of independent and identically distributed random variables. Define $\eta(x) = p(x)f(x)$ and $\hat{\eta}(x) = \frac{1}{n} \sum_{j=1}^n \delta_j W_h(X_j - x)$ as its kernel estimator, where $W_h(u) = h^{-d_x} W(u/h)$. Then,

$$T_n = \frac{1}{n} \sum_{i=1}^n (1 - \delta_i) \frac{(1/n) \sum_{j=1}^n \delta_j W_h(X_j - X_i) \{ g_u(X_i, Y_j, \theta_0) - m_{g_u}(X_j) \}}{\eta(X_i)}$$

$$+ \frac{1}{n} \sum_{i=1}^n (1 - \delta_i) \{ \hat{m}_{g_u}(X_i) - m_{g_u}(X_i) \} \frac{\eta(X_i) - \hat{\eta}(X_i)}{\eta(X_i)}$$

$$+ \frac{1}{n} \sum_{i=1}^n (1 - \delta_i) \left\{ \frac{(1/n) \sum_{j=1}^n \delta_j W_h(X_j - X_i)(m_{g_u}(X_j) - m_{g_u}(X_i))}{\eta(X_i)} \right\}$$

$$:= T_{n1} + T_{n2} + T_{n3}.$$

We now derive the asymptotic distribution of $T_{n1}$. Note that, by exchanging the summation operator,

$$T_{n1} = \frac{1}{n} \sum_{i=1}^n (1 - \delta_i) \frac{(1/n) \sum_{j=1}^n \delta_j W_h(X_j - X_i) \{ g_u(X_i, Y_j, \theta_0) - m_{g_u}(X_j) \}}{\eta(X_i)}$$

$$= \frac{1}{n^2} \sum_{j=1}^{n} \sum_{i=1}^{n} \delta_j \{g_u(X_i, Y_j, \theta) - m_{g_u}(X_j)\} \frac{(1-\delta_i)W_h(X_i - X_j)}{\eta(X_i)}.$$

$$:= \frac{1}{n^2} \sum_{j=1}^{n} \sum_{i=1}^{n} Q_{ij}, \qquad \text{say.}$$

Define

$$\check{T}_{n1} = \frac{1}{n^2} \sum_{j=1}^{n} \sum_{i=1}^{n} E\{Q_{ij} | (X_j, Y_j, \delta_j)\}$$

and write $T_{n1} = \check{T}_{n1} + (T_{n1} - \check{T}_{n1})$. The following derivation will show that $T_{n1}$ is dominated by $\check{T}_{n1}$, while $(T_{n1} - \check{T}_{n1})$ is of smaller order. We note by ignoring terms of $O_p(h^2)$, which are $o_p(n^{-1/2})$ under the assumption that $\sqrt{n}h^2 \to 0$,

$$\check{T}_{n1} = \frac{1}{n} \sum_{j=1}^{n} \delta_j E\left[ \{g_u(X_i, Y_j, \theta_0) - m_{gu}(X_j)\} \frac{(1-\delta_i)W_h(X_i - X_j)}{\eta(X_i)} \Big| X_j, Y_j \right]$$

$$= \frac{1}{n} \sum_{j=1}^{n} \delta_j E_{X_i|X_j, Y_j} \left( E\left[ \{g_u(X_i, Y_j, \theta_0) - m_{gu}(X_j)\} \right. \right.$$

$$\left. \left. \times \frac{(1-\delta_i)W_h(X_i - X_j)}{\eta(X_i)} \Big| X_j, Y_j, X_i \right] \right)$$

$$= \frac{1}{n} \sum_{j=1}^{n} \delta_j E_{X_i|X_j, Y_j} \left[ \{g_u(X_i, Y_j, \theta_0) - m_{gu}(X_j)\} \right.$$

$$\left. \times \frac{(1-P(X_i))W_h(X_i - X_j)}{\eta(X_i)} \right],$$

where $E_{X_i|X_j, Y_j}(\cdot)$ represents conditional expectation on $X_i$ given $(X_j, Y_j)$. Then,

$$\check{T}_{n1} = \frac{1}{n} \sum_{j=1}^{n} \delta_j \int \left[ \{g_u(x, Y_j, \theta_0) - m_{gu}(X_j)\} \frac{\{1-p(x)\}W_h(x - X_j)}{\eta(x)} \right] f(x)\, dx$$

$$= \frac{1}{n} \sum_{j=1}^{n} \delta_j \int \left[ \{g_u(x, Y_j, \theta_0) - m_{gu}(X_j)\} \frac{\{1-p(x)\}}{p(x)} W_h(x - X_j) \right] dx$$

$$= \frac{1}{n} \sum_{j=1}^{n} \delta_j \int \left[ \{g_u(X_j + hs, Y_j, \theta_0) - m_{gu}(X_j)\} \frac{\{1-p(X_j + hs)\}}{p(X_j + hs)} W(s) \right] ds.$$

Since both $g_u$ and $\rho(x) = \{1 - p(x)\}/p(x)$ have bounded second derivative on $x$, and $\sqrt{n}h^2 \to 0$ as $n \to \infty$, a Taylor expansion around $X_j$ leads to

$$\text{(A.1)} \quad \check{T}_{n1} = \frac{1}{n}\sum_{j=1}^{n}\delta_j\{g_u(X_j, Y_j, \theta_0) - m_{gu}(X_j)\}\frac{1 - p(X_j)}{p(X_j)} + o_p(n^{-1/2}).$$

Now we show $T_{n1} - \check{T}_{n1} = o_p(n^{-1/2})$. Let

$$T_{n1i} = \frac{1}{n}\sum_{j=1}^{n}Q_{ij}$$

and

$$\check{T}_{n1i} = \frac{1}{n}\sum_{j=1}^{n}E\{Q_{ij} \mid (X_j, Y_j, \delta_j)\}.$$

By straight forward derivations,

$$\begin{aligned}
&nE(T_{n1} - \check{T}_{n1})^2 \\
\text{(A.2)} \quad &= \frac{1}{n}\sum_{i=1}^{n}E(T_{n1i} - \check{T}_{n1i})^2 + \frac{2}{n}\sum_{i\neq k}E\{(T_{n1i} - \check{T}_{n1i})(T_{n1k} - \check{T}_{n1k})\} \\
&= E(T_{n1i} - \check{T}_{n1i})^2.
\end{aligned}$$

The last step used the fact that $E_{i\neq j}\{(T_{n1i} - \check{T}_{n1i})(T_{n1j} - \check{T}_{n1j})\} = 0$, which can be shown by conditioning on the completely observed portion of data. Thus,

$$\begin{aligned}
&nE(T_{n1} - \check{T}_{n1})^2 \\
&= ET_{n1i}^2 - E\check{T}_{n1i}^2 \\
&\leq ET_{n1i}^2 \\
&\leq E\left\{\frac{(1/n)\sum_{j=1}^{n}\delta_j W_h(X_j - X_i)\{g_u(X_i, Y_j, \theta_0) - m_{gu}(X_j)\}}{\eta(X_i)}\right\}^2 \\
&\to 0,
\end{aligned}$$

by a standard derivation in kernel estimation. This suggests that $T_{n1} = \check{T}_{n1} + o_p(n^{-1/2})$. By a standard argument, it may be shown that $T_{n2} = o_p(n^{-1/2})$. For $T_{n3}$, a similar derivation to that for $T_{n1}$ shows that $T_{n3} = o_p(n^{-1/2})$. Thus,

$$\text{(A.3)} \quad \sqrt{n}T_n \xrightarrow{\mathcal{L}} N[0, E\{(1 - p(X))^2\sigma_{g_u}^2(X)/p(X)\}],$$

where $\sigma_{g_u}^2(X) = \text{Var}\{g_u(X, Y, \theta_0) \mid X\}$.

Note that $\sqrt{n}S_n \xrightarrow{\mathcal{L}} N[0, E\{p(X)\sigma^2_{g_u}(X)\}]$ and $\sqrt{n}R_n \xrightarrow{\mathcal{L}} N[0, \text{Var}\{m_{g_u}(X)\}]$ by the central limit theorem. Furthermore, it can be shown that

$$n \text{Cov}(S_n, T_n) = E\{(1 - p(X))\sigma^2_{g_u}(X)\} + o(1),$$

$n \text{Cov}(R_n, S_n) = 0$ and $n \text{Cov}(R_n, T_n) = o(1)$. It readily follows that

$$\sqrt{n} \begin{pmatrix} S_n \\ T_n \\ R_n \end{pmatrix} \xrightarrow{\mathcal{L}} N\left(0, \begin{bmatrix} \Upsilon & 0 \\ 0 & \text{Var}(m_{g_u}(X)) \end{bmatrix}\right),$$

where

$$\Upsilon = \begin{bmatrix} E\{p(X)\sigma^2_{g_u}(X)\} & E\{(1 - p(X))\sigma^2_{g_u}(X)\} \\ E\{(1 - p(X))\sigma^2_{g_u}(X)\} & E\{(1 - p(X))^2\sigma^2_{g_u}(X)/p(X)\} \end{bmatrix}.$$

Hence, we have

(A.4)    $\sqrt{n}(S_n + T_n + R_n) \xrightarrow{\mathcal{L}} N[0, E\{\sigma^2_{g_u}(X)/p(X)\} + \text{Var}\{m_{g_u}(X)\}].$

Now we consider the asymptotic distribution of

$$A_n = \frac{1}{n}\sum_{i=1}^{n}(1 - \delta_i)\left\{\kappa^{-1}\sum_{\nu=1}^{\kappa} g_u(X_i, \tilde{Y}_{i\nu}, \theta_0) - \hat{m}_{g_u}(X_i)\right\}.$$

Given all the original observations, $n^{-1/2}(1 - \delta_i)\{\kappa^{-1}\sum_{\nu=1}^{\kappa} g_u(X_i, \tilde{Y}_{i\nu}, \theta) - \hat{m}(X_i)\}$, $i = 1, 2, \ldots, n$, are independent with conditional mean zero and conditional variance $(n\kappa)^{-1}(1 - \delta_i)\{\hat{\gamma}_{g_u}(X_i) - \hat{m}^2_{g_u}(X_i)\}$. Here

$$\hat{\gamma}_{g_u}(x) = \sum_{j=1}^{n} \delta_j W_h(x - X_j) g^2_u(x, Y_j, \theta_0)/\hat{\eta}(x)$$

is a kernel estimator of $\gamma_{g_u}(x) = E\{g^2_u(X, Y, \theta_0)|X = x\}$. By verifying Lyapounov's condition, we can show that conditioning on the original observations, $\sqrt{n}A_n$ has an asymptotic normal distribution with mean zero and variance $(n\kappa)^{-1}\sum_{i=1}^{n}(1 - \delta_i)\{\hat{\gamma}_{g_u}(X_i) - \hat{m}^2_{g_u}(X_i)\}$. The conditional variance

(A.5)    $(n\kappa)^{-1}\sum_{i=1}^{n}(1 - \delta_i)\{\hat{\gamma}_{g_u}(X_i) - \hat{m}^2_{g_u}(X_i)\} \xrightarrow{p} \kappa^{-1}E[\{1 - p(X)\}\sigma^2_{g_u}(X)].$

By Proposition A.1, we can show that, as $n \to \infty$ and $\kappa \to \infty$, $\sqrt{n}(S_n + T_n + R_n + A_n)$ converges to a normal distribution with mean 0 and variance

$$\text{Var}\{m_{g_u}(Z, \theta_0)\} + E\{p^{-1}(X)\sigma^2_{g_u}(X)\} = u^\tau \Gamma u.$$

Hence, $n^{-1/2}\sum_{i=1}^{n} \tilde{g}_u(X_i, \theta_0) \xrightarrow{\mathcal{L}} N(0, u^\tau\Gamma u)$. And Lemma A.1 is proved by using the Cramèr–Wold device.  □

LEMMA A.2. *Under the conditions* C1–C5, *as* $n \to \infty$ *and* $\kappa \to \infty$,

$$\frac{1}{n}\sum_{i=1}^{n}\tilde{g}_i(\theta_0)\tilde{g}_i^{\tau}(\theta_0) \xrightarrow{p} \tilde{\Gamma},$$

*where* $\tilde{\Gamma} = E\{p(X)\operatorname{Var}(g|X) + E(g|X)E(g^{\tau}|X)\}$.

PROOF. Consider the $(j,k)$th element of the matrix $\frac{1}{n}\sum_{i=1}^{n}\tilde{g}_i(\theta_0)\tilde{g}_i^{\tau}(\theta_0)$, that is,

$$\frac{1}{n}\sum_{i=1}^{n}\tilde{g}_{i(j)}(\theta_0)\tilde{g}_{i(k)}(\theta_0),$$

where $\tilde{g}_{i(j)}(\theta_0)$ and $\tilde{g}_{i(k)}(\theta_0)$ represent the $j$th and $k$th element of the vector $\tilde{g}_i(\theta_0)$, respectively, for $0 \le j, k \le r$. Similarly, we use $g_{(j)}$ to represent the $j$th element of $g$. Note that

$$\begin{aligned}
\frac{1}{n}&\sum_{i=1}^{n}\tilde{g}_{i(j)}(\theta_0)\tilde{g}_{i(k)}(\theta_0)\\
&=\frac{1}{n}\sum_{i=1}^{n}\delta_i g_{(j)}(Z_i,\theta_0)g_{(k)}(Z_i,\theta_0)\\
&\quad+\frac{1}{n}\sum_{i=1}^{n}(1-\delta_i)\left\{\kappa^{-1}\sum_{\nu=1}^{\kappa}g_{(j)}(X_i,\tilde{Y}_{i\nu},\theta_0)\right\}\left\{\kappa^{-1}\sum_{\nu=1}^{\kappa}g_{(k)}(X_i,\tilde{Y}_{i\nu},\theta_0)\right\}\\
&:=B_{n1}+B_{n2}.
\end{aligned}$$

Moreover,

$$\begin{aligned}
B_{n1} &= \frac{1}{n}\sum_{i=1}^{n}\delta_i\{g_{(j)}(Z_i,\theta_0)-m_{g_{(j)}}(X_i)\}\{g_{(k)}(Z_i,\theta_0)-m_{g_{(k)}}(X_i)\}\\
&\quad-\frac{1}{n}\sum_{i=1}^{n}\delta_i m_{g_{(j)}}(X_i)m_{g_{(k)}}(X_i)+\frac{1}{n}\sum_{i=1}^{n}\delta_i g_{(j)}(Z_i,\theta_0)m_{g_{(k)}}(X_i)\\
&\quad+\frac{1}{n}\sum_{i=1}^{n}\delta_i g_{(k)}(Z_i,\theta_0)m_{g_{(j)}}(X_i)\\
&:=B_{n1a}+B_{n1b}+B_{n1c}+B_{n1d}.
\end{aligned}$$

It is obvious that $B_{n1a}$, $B_{n1b}$, $B_{n1c}$ and $B_{n1d}$ are all sums of independent and identically distributed random variables. By law of large numbers and the continuous mapping theorem, we can show that

$$B_{n1}\xrightarrow{p} E[p(X)\operatorname{Cov}\{g_{(j)}(Z,\theta_0),g_{(k)}(Z,\theta_0)|X\}+p(X)m_{g_{(j)}}(X)m_{g_{(k)}}(X)],$$

where

$$\text{Cov}\{g_{(j)}(Z,\theta_0), g_{(k)}(Z,\theta_0)|X\}$$
$$= E[\{g_{(j)}(Z,\theta_0) - m_{g_{(j)}}(X_i)\}\{g_{(k)}(Z,\theta_0) - m_{g_{(k)}}(X_i)\}|X].$$

Note that

$$B_{n2} = \frac{1}{n}\sum_{i=1}^{n}(1-\delta_i)\left[\left\{\kappa^{-1}\sum_{\nu=1}^{\kappa}g_{(j)}(X_i,\tilde{Y}_{i\nu},\theta_0)\right\}\left\{\kappa^{-1}\sum_{\nu=1}^{\kappa}g_{(k)}(X_i,\tilde{Y}_{i\nu},\theta_0)\right\}\right.$$
$$\left. - \hat{m}_{g_{(j)}}(X_i)\hat{m}_{g_{(k)}}(X_i)\right]$$
$$+ \frac{1}{n}\sum_{i=1}^{n}(1-\delta_i)\{\hat{m}_{g_{(j)}}(X_i)\hat{m}_{g_{(k)}}(X_i) - m_{g_{(j)}}(X_i)m_{g_{(k)}}(X_i)\}$$
$$+ \frac{1}{n}\sum_{i=1}^{n}(1-\delta_i)m_{g_{(j)}}(X_i)m_{g_{(k)}}(X_i)$$
$$:= B_{n2a} + B_{n2b} + B_{n2c}.$$

As $\kappa^{-1}\sum_{\nu=1}^{\kappa}g_{(j)}(X_i,\tilde{Y}_{i\nu},\theta_0)$ has conditional mean $\hat{m}_{g_{(j)}}(X_i)$ given the original observations $\boldsymbol{\chi}_n$, it can be shown that $B_{n2a} \xrightarrow{p} 0$ as $\kappa \to \infty$. By argument similar to those used for (A.3), $B_{n2b} \xrightarrow{p} 0$ as $n \to \infty$. Obviously $\boldsymbol{B}_{n2c}$ is the sum of independent and identically distributed random variables, which leads to $B_{n2c} \xrightarrow{p} E[\{1-p(X)\}m_{g_{(j)}}(X_i)m_{g_{(k)}}(X_i)]$. Hence, we have $B_{n2} \xrightarrow{p} E[\{1-p(X)\}m_{g_{(j)}}(X_i)m_{g_{(k)}}(X_i)]$ as $n \to \infty$ and $\kappa \to \infty$. Therefore,

$$B_{n1} + B_{n2} \xrightarrow{p} E[p(X)\text{Cov}\{g_{(j)}(Z,\theta_0), g_{(k)}(Z,\theta_0)|X\} + m_{g_{(j)}}(X)m_{g_{(k)}}(X)].$$

This completes the proof of Lemma A.2.  □

PROOF OF THEOREM 1.  The proof of Theorem 1 is very similar to that of Lemma 1 of Qin and Lawless [23]. Briefly, we can show

$$t(\theta) = \left\{\frac{1}{n}\sum_{i=1}^{n}\tilde{g}_i(\theta)\tilde{g}_i^{\tau}(\theta)\right\}^{-1}\left\{\frac{1}{n}\sum_{i=1}^{n}\tilde{g}_i(\theta)\right\} + o(n^{-1/3})$$
$$= O(n^{-1/3})$$

almost surely uniformly for all $\theta$ such that $\|\theta - \theta_0\| \le Cn^{-1/3}$ for a positive constant $C$.

From this and Taylor expansion, we can show $\ell_n(\theta) = O(n^{1/3})$ and $\ell_n(\theta_0) = O(\log\log n)$ almost surely. Noting that $\ell(\theta)$ is a continuous function about $\theta$ as $\theta$ belongs to the ball $\|\theta - \theta_0\| \le Cn^{-1/3}$, with probability tending to 1, $\ell_n(\theta)$ has a minimum $\hat{\theta}_n$ in the interior of the ball, and this $\hat{\theta}_n$ satisfies

$$
\begin{aligned}
\frac{\partial \ell_n(\theta)}{\partial \theta}\bigg|_{\theta=\hat{\theta}_n} &= \sum_i \frac{\{\partial t^\tau(\theta)/\partial \theta\}\tilde{g}_i(\theta) + \{\partial \tilde{g}_i(\theta)/\partial \theta\}^\tau t(\theta)}{1 + t^\tau(\theta)\tilde{g}_i(\theta)}\bigg|_{\theta=\hat{\theta}_n} \\
&= \sum_i \frac{1}{1 + t^\tau(\theta)\tilde{g}_i(\theta)}\left\{\frac{\partial \tilde{g}_i(\theta)}{\partial \theta}\right\}^\tau t(\theta)|_{\theta=\hat{\theta}_n} \\
&= 0.
\end{aligned}
$$

Hence, the $\hat{\theta}_n$ satisfies the second equation of (6). From the algorithm of the empirical likelihood formulation, the $\hat{\theta}_n$ automatically satisfies the first equation of (6). This completes the proof of Theorem 1. $\square$

PROOF OF THEOREM 2.   Recall that $\hat{\theta}_n$ and $\hat{t} = t(\hat{\theta}_n)$ satisfy

$$
Q_{1n}(\hat{\theta}_n, \hat{t}) = 0, \qquad Q_{2n}(\hat{\theta}_n, \hat{t}) = 0.
$$

Taking the derivatives with regard to $\theta$ and $t^\tau$,

$$
\begin{aligned}
\frac{\partial Q_{1n}(\theta, 0)}{\partial \theta} &= \frac{1}{n}\sum_i \frac{\partial \tilde{g}_i(\theta)}{\partial \theta}, & \frac{\partial Q_{1n}(\theta, 0)}{\partial t^\tau} &= -\frac{1}{n}\sum_i \tilde{g}_i(\theta)\tilde{g}_i^\tau(\theta), \\
\frac{\partial Q_{2n}(\theta, 0)}{\partial \theta} &= 0, & \frac{\partial Q_{2n}(\theta, 0)}{\partial t^\tau} &= \frac{1}{n}\sum_i \left\{\frac{\partial \tilde{g}_i(\theta)}{\partial \theta}\right\}^\tau.
\end{aligned}
$$

Expanding $Q_{1n}(\hat{\theta}_n, \hat{t})$ and $Q_{2n}(\hat{\theta}_n, \hat{t})$ at $(\theta_0, 0)$, we have

$$
\begin{aligned}
0 &= Q_{1n}(\hat{\theta}_n, \hat{t}) \\
&= Q_{1n}(\theta_0, 0) + \frac{\partial Q_{1n}(\theta_0, 0)}{\partial \theta}(\hat{\theta}_n - \theta_0) + \frac{\partial Q_{1n}(\theta_0, 0)}{\partial t^\tau}(\hat{t} - 0) + o_p(\zeta_n), \\
0 &= Q_{2n}(\hat{\theta}_n, \hat{t}) \\
&= Q_{2n}(\theta_0, 0) + \frac{\partial Q_{2n}(\theta_0, 0)}{\partial \theta}(\hat{\theta}_n - \theta_0) + \frac{\partial Q_{2n}(\theta_0, 0)}{\partial t^\tau}(\hat{t} - 0) + o_p(\zeta_n),
\end{aligned}
$$

where $\zeta_n = \|\hat{\theta}_n - \theta_0\| + \|\hat{t}\|$. Then, we have

$$
\begin{pmatrix} \hat{t} \\ \hat{\theta}_n - \theta_0 \end{pmatrix} = S_n^{-1}\begin{pmatrix} -Q_{1n}(\theta_0, 0) + o_p(\zeta_n) \\ o_p(\zeta_n) \end{pmatrix},
$$

where

$$
S_n = \begin{pmatrix} \dfrac{\partial Q_{1n}}{\partial t^\tau} & \dfrac{\partial Q_{1n}}{\partial \theta} \\ \dfrac{\partial Q_{2n}}{\partial t^\tau} & 0 \end{pmatrix}_{(\theta_0, 0)} \xrightarrow{p} \begin{pmatrix} S_{11} & S_{12} \\ S_{21} & 0 \end{pmatrix} = \begin{pmatrix} -\tilde{\Gamma} & E\left(\dfrac{\partial g}{\partial \theta}\right) \\ E\left(\dfrac{\partial g}{\partial \theta}\right)^\tau & 0 \end{pmatrix}.
$$

Here $\partial Q_{1n}/\partial t^{\tau}|_{(\theta_0,0)} \overset{p}{\to} S_{11}$ follows from Lemma A.2, and $\partial Q_{1n}/\partial \theta|_{(\theta_0,0)} \overset{p}{\to} S_{12}$ can be derived by arguments similar to those used for the proof of Lemma A.1. Note that $Q_{1n}(\theta_0,0) = \frac{1}{n}\sum_{i=1}^{n}\tilde{g}_i(\theta_0) = O_p(n^{-1/2})$, it follows that $\zeta_n = O_p(n^{-1/2})$. After some matrix manipulation, we have

$$\sqrt{n}(\hat{\theta}_n - \theta_0) = S_{22.1}^{-1}S_{21}S_{11}^{-1}\sqrt{n}Q_{1n}(\theta_0,0) + o_p(1),$$

where $V = S_{22.1}^{-1} = \{E(\frac{\partial g}{\partial \theta})^{\tau}\tilde{\Gamma}^{-1}E(\frac{\partial g}{\partial \theta})\}^{-1}$. By Lemma A.1, $\sqrt{n}Q_{1n}(\theta_0,0) \overset{\mathcal{L}}{\to} N(0,\Gamma)$, and the theorem follows.  $\square$

PROOF OF THEOREM 3.   Notice that

$$\mathcal{R}(\theta_0) = 2\left[\sum_i \log\{1 + t_0^{\tau}\tilde{g}_i(\theta_0)\} - \sum_i \log\{1 + \hat{t}^{\tau}\tilde{g}_i(\hat{\theta}_n)\}\right],$$

where $t_0 = t(\theta_0)$, and

$$\ell(\hat{\theta}_n, \hat{t}) = \sum_i \log\{1 + \hat{t}^{\tau}\tilde{g}_i(\hat{\theta}_n)\} = -\frac{n}{2}Q_{1n}^{\tau}(\theta_0,0)AQ_{1n}(\theta_0,0) + o_p(1),$$

where $A = S_{11}^{-1}(I + S_{12}S_{22.1}^{-1}S_{21}S_{11}^{-1})$. Under $H_0$,

$$\frac{1}{n}\sum_i \frac{1}{1 + t_0^{\tau}\tilde{g}_i(\theta_0)}\tilde{g}_i(\theta_0) = 0, \qquad t_0 = -S_{11}^{-1}Q_{1n}(\theta_0,0)S_{11}^{-1}Q_{1n}(\theta_0,0) + o_p(1)$$

and $\sum_i \log\{1 + t_0^{\tau}\tilde{g}_i(\theta_0)\} = -\frac{n}{2}Q_{1n}^{\tau}(\theta_0,0)S_{11}^{-1}Q_{1n}(\theta_0,0) + o_p(1)$. Thus,

$$\begin{aligned}
\mathcal{R}(\theta_0) &= nQ_{1n}^{\tau}(\theta_0,0)(A - S_{11}^{-1})Q_{1n}(\theta_0,0) + o_p(1) \\
&= \sqrt{n}Q_{1n}^{\tau}(\theta_0,0)S_{11}^{-1}S_{12}S_{22.1}^{-1}S_{21}S_{11}^{-1}\sqrt{n}Q_{1n}(\theta_0,0) + o_p(1).
\end{aligned}$$

Note that

$$S_{11}^{-1}S_{12}S_{22.1}^{-1}S_{21}S_{11}^{-1} \overset{p}{\to} \tilde{\Gamma}^{-1}E\left(\frac{\partial g}{\partial \theta}\right)VE\left(\frac{\partial g}{\partial \theta}\right)^{\tau}\tilde{\Gamma}^{-1}$$

and, by Lemma A.1, $\sqrt{n}Q_{1n}(\theta_0,0) \overset{\mathcal{L}}{\to} N(0,\Gamma)$. This implies the theorem.  $\square$

PROOF OF THEOREM 4.   The proof of Theorem 4 essentially involves establishing the bootstrap version of Lemma A.1 to Theorem 3. We only outline the main steps in proving the bootstrap version of Lemma A.1 here.

Let $X_i^*$, $Y_i^*$, $\tilde{Y}_{iv}^*$, $\delta_i^*$, $\tilde{g}_{ui}^*$ be counter parts of $X_i$, $Y_i$, $\tilde{Y}_{iv}$, $\delta_i$, $\tilde{g}_{ui}$ in the bootstrap sample; and $S_n(\hat{\theta}_n)$, $A_n(\hat{\theta}_n)$, $T_n(\hat{\theta}_n)$ and $R_n(\hat{\theta}_n)$ be the quantities $S_n$, $A_n$, $T_n$ and $R_n$ with $\theta_0$ replaced by $\hat{\theta}_n$, respectively. Furthermore, let $S_n^*(\hat{\theta}_n)$, $A_n^*(\hat{\theta}_n)$, $T_n^*(\hat{\theta}_n)$ and $R_n^*(\hat{\theta}_n)$ be their bootstrap counterparts. First, we will consider the conditional distribution of $\sqrt{n}\{S_n^*(\hat{\theta}_n) + T_n^*(\hat{\theta}_n) + R_n^*(\hat{\theta}_n) - S_n(\hat{\theta}_n) - T_n(\hat{\theta}_n) - R_n(\hat{\theta}_n)\}$ given the original data. We use $E_*(\cdot)$ and $\mathrm{Var}_*(\cdot)$ to represent the conditional expectation and variance given the original data, respectively. Define

$$\hat{m}_{gu}(x, \hat{\theta}_n) = \frac{\sum_{i=1}^n \delta_i W((x - X_i)/h) g_u(x, Y_i, \hat{\theta}_n)}{\sum_{i=1}^n \delta_i W((x - X_i)/h)}$$

and

$$\hat{m}_{gu}^*(x, \hat{\theta}_n) = \frac{\sum_{i=1}^n \delta_i^* W((x - X_i^*)/h) g_u(x, Y_i^*, \hat{\theta}_n)}{\sum_{i=1}^n \delta_i^* W((x - X_i^*)/h)}.$$

Note that $S_n(\hat{\theta}_n) + T_n(\hat{\theta}_n) + R_n(\hat{\theta}_n) = \frac{1}{n} \sum_{i=1}^n \{\delta_i g_u(Z_i, \hat{\theta}_n) + (1 - \delta_i) \hat{m}_{gu}(X_i, \hat{\theta}_n)\}$. Thus,

$$S_n^*(\hat{\theta}_n) + T_n^*(\hat{\theta}_n) + R_n^*(\hat{\theta}_n) - S_n(\hat{\theta}_n) - T_n(\hat{\theta}_n) - R_n(\hat{\theta}_n)$$

$$= \frac{1}{n} \sum_{i=1}^n [\delta_i^* g_u(Z_i^*, \hat{\theta}_n) + (1 - \delta_i^*) \hat{m}_{gu}(X_i^*, \hat{\theta}_n)$$

$$- E_* \{\delta_i^* g_u(Z_i^*, \hat{\theta}_n) + (1 - \delta_i^*) \hat{m}_{gu}(X_i^*, \hat{\theta}_n)\}]$$

$$+ \frac{1}{n} \sum_{i=1}^n (1 - \delta_i^*) \{\hat{m}_{gu}^*(X_i^*, \hat{\theta}_n) - \hat{m}_{gu}(X_i^*, \hat{\theta}_n)\}$$

$$:= B_1 + B_2.$$

It can be shown that $B_2 = o_p(n^{-1/2})$. For $B_1$, we can apply the central limit theorem for bootstrap samples, for example, [31] to show that the conditional distribution of $\sqrt{n} B_1$ given $\boldsymbol{\chi}_n$ is asymptotically normal with mean zero and variance $\mathrm{Var}_* \{\delta_i^* g_u(Z_i^*, \hat{\theta}_n) + (1 - \delta_i^*) \hat{m}_{gu}(X_i^*, \hat{\theta}_n)\}$.

Using similar methods for Lemma A.1, we can also derive the conditional distribution of $\sqrt{n} A_n^*(\hat{\theta}_n)$ given the observations in the bootstrap sample that are not imputed. Then by employing Proposition A.1, it follows that the conditional distribution of $n^{-1/2} \sum_{i=1}^n \tilde{g}_{ui}^*(\hat{\theta}_n)$ given $\boldsymbol{\chi}_n$ is asymptotically normal with mean zero and variance $\hat{\sigma}_u^{2*} = \mathrm{Var}_* \{\delta_i^* g_u(Z_i^*, \hat{\theta}_n) + (1 - \delta_i^*) \hat{m}_{gu}(X_i^*, \hat{\theta}_n)\}$. The bootstrap version of Lemma A.1 is justified by noting that $\hat{\sigma}_u^{2*}$ converges in probability to $u^\tau \Gamma u$ as $n \to \infty$, then employing the Cramèr–Wold device. $\square$

**Acknowledgments.** We thank two referees for constructive comments and suggestions which improve the presentation of the paper. The paper is part of Ph.D. thesis of Dong Wang at Department of Statistics, Iowa State University. He thanks Dan Nettleton for support.

# REFERENCES

[1] Azzalini, A. (1981). A note on the estimation of a distribution function and quantiles by a kernel method. *Biometrika* **68** 326–328. MR0614972

[2] Azzalini, A. and Capitanio, A. (2003). Distributions generated by perturbation of symmetry with emphasis on a multivariate skew *t*-distribution. *J. Roy. Statist. Soc. Ser. B* **65** 367–389. MR1983753

[3] Boos, D. D. (1992). On generalized score tests. *Amer. Statist.* **46** 327–333.

[4] Bowman, A., Hall, P. and Prvan, T. (1998). Bandwidth selection for the smoothing of distribution functions. *Biometrika* **85** 799–808. MR1666695

[5] Chen, S. X. and Cui, H. J. (2007). On the second-order properties of empirical likelihood with moment restrictions. *J. Econometrics* **141** 492–516.

[6] Chen, X., Hong, H. and Tarozzi, A. (2008). Semiparametric efficiency in GMM models with auxiliary data. *Ann. Statist.* **36** 808–843. MR2396816

[7] Cheng, P. E. (1994). Nonparametric estimation of mean functionals with data missing at random. *J. Amer. Statist. Assoc.* **89** 81–87.

[8] Cramér, H. (1946). *Mathematical Methods of Statistics.* Princeton Univ. Press, Princeton, NJ. MR0016588

[9] Godambe, V. P. (1991). *Estimating Functions.* Oxford Univ. Press, Oxford. MR1163992

[10] Hall, P. and Murison, R. D. (1993). Correcting the negativity of high-order kernel density estimators. *J. Multivariate Anal.* **47** 103–122. MR1239108

[11] Hansen, L. (1982). Large sample properties of generalized method of moment estimators. *Econometrica* **50** 1029–1084. MR0666123

[12] Härdle, W. (1990). *Applied Nonparametric Regression.* Cambridge Univ. Press, Cambridge. MR1161622

[13] Hjort, N. L., McKeague, I. and Van Keilegom, I. (2008). Extending the scope of empirical likelihood. *Ann. Statist.* To appear.

[14] Ibrahim, J. G., Chen, M., Lipsitz, S. R. and Herring, A. H. (2005). Missing-data methods for generalized linear models: A comparative review. *J. Amer. Statist. Assoc.* **100** 332–346. MR2166072

[15] Jin, Z. and Shao, Y. (1999). On kernel estimation of a multivariate distribution function. *Statist. Probab. Lett.* **2** 163–168. MR1665267

[16] Kolaczyk, E. D. (1994). Empirical likelihood for generalized linear models. *Statist. Sinica* **4** 199–218. MR1282871

[17] Little, R. J. A. and Rubin, D. B. (2002). *Statistical Analysis with Missing Data*, 2nd ed. Wiley, Hoboken, NJ. MR1925014

[18] McCullagh, P. and Nelder, J. A. (1983). *Generalized Linear Models.* Chapman and Hall, London. MR0727836

[19] Owen, A. (1988). Empirical likelihood ratio confidence intervals for a single functional. *Biometrika* **75** 237–249. MR0946049

[20] Owen, A. (1990). Empirical likelihood ratio confidence regions. *Ann. Statist.* **18** 90–120. MR1041387

[21] Owen, A. B. (2001). *Empirical Likelihood.* Chapman and Hall/CRC, Boca Raton, FL.

[22] Pierce, J. L., Lu, L., Gu, J., Silver, L. M. and Williams, R. W. (2004). A new set of BXD recombinant inbred lines from advanced intercross populations in mice. *BMG Genetics* **5** http://www.biomedcentral.com/1471-2156/5/7.

[23] Qin, J. and Lawless, J. (1994). Empirical likelihood and general estimating equations. *Ann. Statist.* **22** 300–325. MR1272085

[24] ROBINS, J. M. and ROTNITZKY, A. (1995). Semiparametric efficiency in multivariate regression models with missing data. *J. Amer. Statist. Assoc.* **90** 122–129. MR1325119

[25] ROBINS, J. M., ROTNITZKY, A. and ZHAO, L. P. (1994). Estimation of regression coefficients when some regressors are not always observed. *J. Amer. Statist. Assoc.* **89** 846–866. MR1294730

[26] ROBINS, J. M., ROTNITZKY, A. and ZHAO, L. P. (1995). Analysis of semiparametric regression models for repeated outcomes in the presence of missing data. *J. Amer. Statist. Assoc.* **90** 106–121. MR1325118

[27] ROSENBAUM, P. R. and RUBIN, D. B. (1983). The central role of the propensity score in observational studies for causal effects. *Biometrika* **70** 41–55. MR0742974

[28] RUBIN, D. B. (1987). *Multiple Imputation for Nonresponse in Surveys.* Wiley, New York. MR0899519

[29] SCHENKER, N. and WELSH, A. H. (1988). Asymptotic results for multiple imputation. *Ann. Statist.* **16** 1550–1566. MR0964938

[30] SHAO, J. and SITTER, R. R. (1996). Bootstrap for imputed survey data. *J. Amer. Statist. Assoc.* **91** 1278–1288. MR1424624

[31] SHAO, J. and TU, D. (1985). *The Jackknife and Bootstrap.* Springer, New York. MR1351010

[32] TSIATIS, A. A. (2006). *Semiparametric Theory and Missing Data.* Springer, New York. MR2233926

[33] WANG, Q. and RAO, J. N. K. (2002). Empirical likelihood-based inference under imputation for missing response data. *Ann. Statist.* **30** 896–924. MR1922545

[34] WANG, J., WILLIAMS, R. W. and MANLY, K. F. (2003). WebQTL: Web-based complex trait analysis. *Neuroinformatics* **1** 299–308.

[35] WILLIAMS, R. W., GU, J., QI, S. and LU, L. (2001). The genetic structure of recombinant inbred mice: High-resolution consensus maps for complex trait analysis. *Genome Biology* **2** RESEARCH0046.

[36] YATES, F. (1933). The analysis of replicated experiments when the field results are incomplete. *Emporium Journal of Experimental Agriculture* **1** 129–142.

DEPARTMENT OF STATISTICS
UNIVERSITY OF NEBRASKA-LINCOLN
LINCOLN, NEBRASKA 68583
USA
E-MAIL: dwang3@unl.edu

DEPARTMENT OF STATISTICS
IOWA STATE UNIVERSITY
AMES, IOWA 50011
USA
AND
DEPARTMENT OF BUSINESS
STATISTICS AND ECONOMICS
GUANGHUA SCHOOL OF MANAGEMENT
PEKING UNIVERSITY
PEKING 100871
CHINA
E-MAIL: songchen@iastate.edu